\documentclass[11pt,a4paper]{article}
\usepackage{wrapfig,epsfig,latexsym,amssymb}

\newtheorem{Proposition}{\sc Proposition}

\long\def\symbolfootnote[#1]#2{\begingroup\def\thefootnote{\fnsymbol{footnote}}
\footnote[#1]{#2}\endgroup}

\begin{document}

\thispagestyle{empty}

\title{Cartan--Helgason theorem, Poisson transform, and
Furstenberg--Satake compactifications}

\author{Adam Kor\'{a}nyi*}

\date{}

\maketitle

\symbolfootnote[0]{*Partially supported by a PSC--CUNY grant.}

\abstract{The connections between the objects mentioned in the
title are used to give a short proof of the Cartan--Helgason
theorem and a natural construction of the compactifications.}

\section*{Introduction}

For a real semisimple Lie group $G$ we write, as usual, $K$ for a
maximal compact subgroup and $MAN$ for a minimal parabolic. The
Cartan--Helgason theorem consists of two parts, one of these
(Prop. 2 here) states that a finite-dimensional irreducible
representation $(\rho, V)$ of $G$ has a $K$-fixed vector $e$ if
and only if it has an $MN$-fixed vector (which is then a highest
weight vector $v^+$ of $\rho$). Of course, when $\rho $ is a
faithful representation, the orbit $\rho(G)e$ gives an imbedding
of the symmetric space $X\simeq G/K$ into the space $PV$ of lines in
$V$ while $\rho(G)v^+$ gives an imbedding of a space $G/B$ where
$B\supset MAN$ (i.e. $B$ is a parabolic subgroup, and $G/B$ is one
of the Poisson boundaries of $X$ in the sense of [2], [12]). In $PV$
now $G/B$ appears as part of the topological boundary of the image
of $X$.

Two things follow from these observations. First, they give a
natural approach to the Cartan--Helgason theorem [3], [5, p. 535],
[6, p. 139] which we are splitting into Propositions 1 and 2: the
proof of Proposition 2 is based on the Poisson transform. Second,
we can consider the full closure of $X$ in $PV$. It turns out that
this is always one of the compactifications constructed originally
by Satake [13] and reconstructed later in other ways in [2], [12],
[1], [9]. In fact, the construction we sketch here may be
conceptually the simplest of all. As J. A. Wolf tells me, the  idea
of using spherical representations to construct compactifications
had also been suggested by R. Hermann some time ago.
\vskip 1em

1. In the following $\mathfrak{g}$ will be  a real semisimple Lie
algebra, $\mathfrak{g}=\mathfrak{k}+\mathfrak{p}$ its
decomposition under a Cartan involution $\theta, \mathfrak{g}^\mathbb C$
its complexification, $\mathfrak{u}=\mathfrak{k}+{\rm
i}\mathfrak{p}$. We choose a maximal subalgebra
$\mathfrak{a}\subset \mathfrak{p}$ and complete it to a Cartan
subalgebra $\mathfrak{h}=\mathfrak{a}+\mathfrak{t}$ of
$\mathfrak{g}$ (so $\mathfrak{t}\subset \mathfrak{k})$. We
identify $\mathfrak{h}^\mathbb C$ with its dual under the Killing form.
The roots with respect to $\mathfrak{h}^\mathbb C$ span the real form
$\mathfrak{h}_0= \mathfrak{a}+{\rm i}\mathfrak{t}$ of
$\mathfrak{h}^\mathbb C$. The restriction of the Killing form to
$\mathfrak{h}_0$ is positive definite, we denote it by $(. \mid .)$.
Given our identification, a restricted root ($\mathfrak{a}$-root)
is the same as the orthogonal projection of an
$\mathfrak{h}^\mathbb C$-root onto $\mathfrak{a}$. For an
$\mathfrak{h}^\mathbb C$-root $\alpha$ we denote the corresponding root
space in $\mathfrak{g}^\mathbb C$ by $\mathfrak{g}_\alpha$. For an
$\mathfrak{a}$-root $\gamma$ we denote the corresponding root
space by $\mathfrak{g}^\gamma$. We choose an ordering of
$\mathfrak{h}_0$ and we set
$\mathfrak{n}=\sum_{\gamma>0}\mathfrak{g}^\gamma,
\bar{\mathfrak{n}}=\theta \mathfrak{n}$. $G^\mathbb C$ will be the simply
connected group with Lie algebra $\mathfrak{g}^\mathbb C$. The analytic
subgroups of $G^\mathbb C$ for $\mathfrak{g}, \mathfrak{k}, \mathfrak{a},
\mathfrak{n}, \bar{\mathfrak{n}}$ will be denoted $G, K, A, N,
\bar{N}$, while $M, M'$ with Lie algebra $\mathfrak{m}$ will be
the centralizer resp. normalizer of $\mathfrak{a}$ in $K$.
$W=M'/M$ is the Weyl group, $\mathfrak{a}^+$ the open positive
Weyl chamber.

The weights of a finite dimensional representation $(\rho, V)$ of
$\mathfrak{g}^\mathbb C$ (or, what is the same, of $G^\mathbb C$) are in
$\mathfrak{h}_0$. If $\Lambda$ is the highest weight we denote by
$v^+$ a highest weight vector. We will always equip $V$ with a
Hermitian inner product such that $\rho(U)$ is unitary (hence,
$\rho(\mathfrak{a})$ is Hermitian).

The following proposition is the first half of the
Cartan--Helgason theorem. Without claiming any originality, for
reference in Sec. 2, we give a proof based on some fundamental
facts about the structure of $G^\mathbb C$.

\begin{Proposition}  (a) $v^+$ is fixed under the connected
component $M_0$ of $M$ iff $\Lambda \in \mathfrak{a}$.\\ (b) $v^+$
is fixed under $M$ iff, in addition,
\begin{equation}
\frac{(\Lambda | \gamma)}{(\gamma | \gamma)}\in \mathbb{Z}
\end{equation}
for all restricted roots $\gamma$.\\
(c) Any element $\lambda$ of
$\bar{\mathfrak{a}}^+$ satisfying (1) is the highest weight of a
representation of $\mathfrak{g}^\mathbb C$.
\end{Proposition}

\emph{Proof.} (a) $v^+$ is $M_0$-fixed iff $\rho(H)v^+=0\ (\forall H\in {\rm
i}\mathfrak{t})$ and $\rho(X_\alpha)v^+=0$ for all $X_\alpha\in
\mathfrak{g}_\alpha$ such that $\mathfrak{g}_\alpha \subset
\mathfrak{m}$. Now $\rho(H)v^+= (\Lambda| H)v^+=0\ (H\in {\rm
i}\mathfrak{t})$ by itself amounts to $\Lambda\in \mathfrak{a}$.
But $\Lambda\in \mathfrak{a}$ automatically also implies
$\rho(X_\alpha )v^+=0$ for $\mathfrak{g}_\alpha \subset
\mathfrak{m}$, i.e. for $\alpha\perp \Lambda$: In fact by the
weight-string property (e.g. [7, p. 114]) $\Lambda \pm \alpha$ are
either both weights or neither one is, and the first possibility
is excluded by the maximality of $\Lambda$.

(b) By a result of Satake [13] (cf. also [4, p. 435]) $M=Z_1M_0$,
where $Z_1=\exp({\rm i}\mathfrak{a})\cap K$. So we have to show
only that $v^+$ is $Z_1$-fixed iff (1) holds. It is well known
(e.g. [4, p. 322]) that, since $G^\mathbb C$ and therefore $U$ are simply
connected, $\exp {\rm i}H \in K$ for $H$ in $\mathfrak{a}$ is
equivalent to $H$ being in the lattice generated by the vectors
$\frac{\pi}{(\gamma | \gamma)} \gamma$ with the simple restricted
roots $\gamma$. Since $\rho (\exp{\rm i}H)v^+={\rm e}^{{\rm
i(\Lambda \mid H)}}v^+$, $Z_1$-invariance of $v^+$ amounts exactly to
(1).

(c) We only have to check the standard integrality condition with
respect to $\mathfrak{h}^\mathbb C$-roots. For such a root $\alpha$ we
denote by $\bar{\alpha}$ its restriction (projection) to
$\mathfrak{a}$. It is well known [4, p. 322] that
$\frac{(\alpha|\alpha)}{(\bar{\alpha}|\bar{\alpha})}=1, 2$ or 4,
with 4 only when $2\bar{\alpha}$ is also a restricted root. In the
first two cases $2\frac{(\Lambda|\alpha)}{(\alpha|\alpha)}\in
\mathbb{Z}$ trivially. In the third case
$2\frac{(\Lambda|\alpha)}{(\alpha|\alpha)}=\frac{(\Lambda|2\bar{\alpha})}{(2\bar{\alpha}|2\bar{\alpha})}$
is again in $\mathbb{Z}$, by (1). \hfill $\Box$
\vskip 1em

Now we come to the second half of the Cartan--Helgason theorem.
The proof given here is the main point of this article.

\begin{Proposition} An
irreducible representation $(\rho, V)$ of $G$ has a $K$-fixed
vector iff it has an MN-fixed
vector.
\end{Proposition}
\emph{Proof.} Again, the existence of an $MN$-fixed vector amounts
to $v^+ $ being $M$-fixed.

Suppose $v^+$ is $M$-fixed. Then $f(g)=\rho(g)v^+$ is a $V$-valued
function transforming as $f(gm(\exp H)n)= e^{(\Lambda|H)}f(g)$ for
$g\in G, m\in M, H\in \mathfrak{a}, n\in N$. (It is a section
lifted to $G$ of a homogeneous line bundle tensored with $V$.) Its
Poisson transform is
$$(P_{\Lambda+\rho}f)(g)=\int_K\rho(gk)v^+ {\rm d} k $$
(in standard notation, with $\rho$ the half-sum of positive
$\mathfrak{a}$-roots; cf. [14, p. 81]). This can also be written
$\rho(g)e$, where $e=\int_K \rho(k)v^+ {\rm d}k$ is $K$-fixed, we
have to show only that $e\neq 0$. Now the Fatou-type theorem of
Michelson [11] (see also [14, p. 83], [6, p. 120]) says that, for
$H\in \mathfrak{a}^+$.
$$ \lim_{t\rightarrow \infty} e^{-t(\Lambda|H)}\rho (\exp
tH)e= c f(e)=cv^+$$ with $c>0$. It follows that $e\neq 0$.

(We also note that, since $\Lambda\in \mathfrak{a}^+$, the
standard proof of the Fatou-type theorem can be considerably simplified:
The convergence of the integral defining $c$ is obvious without even
using the explicit expression of the Jacobian of the map $\bar{n}
\rightarrow k(\bar{n})$.)

Conversely, suppose there exists a $K$-invariant $e\neq 0$. The
weight spaces $V_\lambda$ for different weights are mutually
orthogonal. We write $e=\sum e_\lambda$ with $e_\lambda\in
V_\lambda$. Now $e_\Lambda\neq 0$, because otherwise $(e \mid v^+)=0$,
hence $(e \mid \rho(k)\rho(a)\rho(n)v^+)=(\rho(k^{-1})e \mid \rho(a)
\rho(n)v^+)=0$ for all $k\in K, a\in A, n\in N$, which is
impossible.

We have $\rho(\exp tH)e=\sum_\lambda e^{t\lambda(H)}e_\lambda$.
For $H\in \mathfrak{a}^+$, we have $\Lambda(H)>\lambda(H)$ for all weights
$\lambda\neq \Lambda$. Hence
$$\lim_{t\rightarrow \infty}e^{-t\Lambda(H)}\rho (\exp tH)e=e_\Lambda.$$
Since $e_\Lambda$ is a limit of $M$-fixed vectors, it is
$M$-fixed.\hfill $\Box$
\vskip 1em 
2. We continue with the setup of the preceding section, we
consider a faithful representation $(\rho, V)$ of $\mathfrak{g}^\mathbb C$
with highest weight $\Lambda\in \mathfrak{a}$. 
When $e\neq 0$ is a $K$-fixed vector,
the map $g\cdot o\longmapsto \rho(g)e$ (we distinguish between the
coset space $G/K$ and the symmetric space $X$; we denote by $o\in
X$ the point corresponding to $K$) is an equivariant imbedding of
$X$ into $V$. This is clear since in each simple factor the
$K$-part is a maximal subgroup. We
write $\tilde{v}$ for the image of $v$ in the projective space $PV$.
Then we also have that $g\cdot o\longmapsto \rho(g)\tilde{e}$ is
an imbedding $X\rightarrow PV$. This is so because $\rho(g)$
is scalar only when $g\in Z$, the center of $G$, and $Z$ is
contained in $K$.

Since we also have $Z\subset M$, we see that $g\longmapsto
(\rho(g)v^+)\tilde{}~$ is an equivariant map $G\rightarrow PV$. The
stabilizer $B$ of $(v^+)\tilde{}~$ contains $MAN$ by Prop. 2 (so is a
parabolic group). So we have $X$ and $G/B$ imbedded in $PV$ and
the proof of Prop. 2 shows that $G/B$ is a part of the topological
boundary of $X$.

We may also note that, denoting by $V'$ the orthocomplement of $e$
in $V$ and imbedding it into $PV$ by $v'\longmapsto (e+v')\tilde{}~$ the
images of $X, G/B$ are actually contained in a bounded part of the
vector space $V'$. This is clear for $\rho(A)\tilde{e}$, since 
$\rho(a)e$ is a positive combination of the orthogonal system formed
by the $e_\lambda$ (cf. the proof of Prop. 2). Then it is also 
true for $\rho(G)\tilde{e}=\rho(K)\rho(A)\tilde{e}$ since $\rho(K)$ 
acts on $V'$ by rotations.

The closure of the image of $X$ in $PV$ (or $V'$) is a
compactification to which the action of $G$ extends naturally. As
we will now indicate, what we get in this way are exactly the
Furstenberg--Satake compactifications [13], [2], [12], [1], [10].
Satake's construction is a special instance of ours, he works with
a special class of representations which he obtains from an
arbitrary representation $\sigma$ as the Cartan product of
$\sigma$ and the contragradient $\sigma^\wedge$ composed with $\theta$,
and which are realized on vector spaces of Hermitian matrices. So
our construction is close in spirit to Satake's.

We denote by $\Pi$ the set of positive restricted roots. $\Lambda$
determines a subset $E_0 \subset \Pi$ defined as those $\gamma\in
\Pi$ which are orthogonal to $\Lambda$. We denote by
$\mathfrak{a}(E_0)$ the common zero-space of the elements of $E_0$ and
 by $\mathfrak{a}(E_0)^+$ the subset where
all $\gamma$ in $\Pi-E_0$ take positve values. So
$\mathfrak{a}(E_0)^+$ is the largest (open) face of
$\overline{\mathfrak{a}^+}$ which is perpendicular to $\Lambda$.

Now we can make Proposition 1 a little more precise. The $\theta$-image
of the subalgebra $\mathfrak{n}^{E_0}=\sum_{\gamma\perp
\mathfrak{a}(E_0)} \mathfrak{g}^\gamma$, and $\mathfrak{a}^{E_0}$,
the orthocomplement of $\mathfrak{a}(E_0)$ in $\mathfrak{a}$,
annihilate $v^+$. Together with $\mathfrak{m}+\mathfrak{n}$ they
form the subalgebra $\mathfrak{m}(E_0)+ \mathfrak{n}(E_0)$
annihilating $v^+$; here $\mathfrak{m}(E_0)= \mathfrak{m}_K(E_0)+
\mathfrak{a}^{E_0}+\mathfrak{n}^{E_0}+\theta \mathfrak{n}^{E_0}$
with $\mathfrak{m}_K(E_0)$ the centralizer of $\mathfrak{a}(E_0)$
in $\mathfrak{k}$ and $\mathfrak{n}(E_0)=\sum_{\gamma\notin
\mathfrak{a}^{E_0}}\mathfrak{g}^\gamma$.

We write $B(E_0)$ for the group generated by the analytic subgroup
corresponding to  $\mathfrak{b}(E_0)= \mathfrak{m}(E_0)+
\mathfrak{a}(E_0)+ \mathfrak{n}(E_0)$ and by $M$. This is the
stabilizer of $(v^+)\tilde{}~$ in $PV$. (As it is well known and easy to
prove, with different choices of $E_0$ these are the only closed
subgroups of $G$ containing $MAN$. The parabolic subgroups are, by
definition, their conjugates.)

To describe the closure of $\rho(G)\tilde{e}$ in $PV$, let now $E$
be any subset of $\Pi$. The orbit $X^E=M(E)\cdot o$ is a symmetric
subspace of $X$. The
imbedding of $X$ in $PV$ induces an imbedding of $X^E$ as
$\rho(M(E))\tilde{e}$. Any point in it can be written as
$\rho(k^E)\rho(a^E)\tilde{e}$ with $k^E\in M_K(E),\ a^E\in A^E$.

We choose an $H\in \mathfrak{a}(E)^+$. We write $e=\sum_\lambda
e_\lambda$ where $\lambda$ runs through the restricted weights of
$\rho$ (unlike in Sec. 1 where we worked with the
$\mathfrak{h}$-weights). Since $\Lambda$ is its own restriction,
still $e_\Lambda\neq 0$. We have, for all $m_E\in M(E)$,
\begin{equation}
\rho(\exp tH)\rho(m_E)e=\rho(m_E)\sum_\lambda
e^{t(\lambda \mid H)} e_\lambda
\end{equation}
As
$t\rightarrow \infty$, the dominating terms in the sum are those
with $\lambda \equiv \Lambda({\rm mod}\ \mathfrak{a}^E)$. The limit, 
on the boundary of the image in $PV$, will therefore be
$\rho(m_E)(e^E)\tilde{}~$, where $e^E=\sum_{\lambda \equiv \Lambda
({\rm mod}\ \mathfrak{a}^E)} e_\lambda$. The limit of the family of
sets $\rho(exp tH) \rho(M(E))\tilde{e}$ will be $\rho(m(E))(e^E)\tilde{}~$.

The restricted weights such that
 $\lambda \equiv \Lambda ({\rm mod}\ \mathfrak{a}^E)$ are those 
that arise in the form
$\Lambda-\sum m_j\delta_j$ with $\delta_j\in E$. It is not hard to
see (cf. [13 , Lemma 8]) that on $V^E$, the direct sum of the
corresponding spaces $V_\lambda$, the restriction of $\rho$ to
$M(E)$ is an irreducible representation, to be denoted $\rho^E$.
Its highest restricted weight (with respect to $\mathfrak{a}^E$)
is the projection $\Lambda^E$ of $\Lambda$ onto $\mathfrak{a}^E$, 
and $e^E$ is an $M_K(E)$-fixed vector of it.

In general $\rho(M(E))(e^E)\tilde{}~$ is not a one-to-one image of $X^E$;
there are in general several subsets $E$ for which $X^E$ has the 
same image. For any $E \subset \Pi$ we have that $X^E$ is the direct 
product of irreducible symmetric spaces $X^{E_i}$. The $E_i$ are called 
the components of $E$, and $E$ is said to be $E_o$-connected if none of 
its components is entirely contained in $E_o$. Clearly, it is exactly
when $E$ is $E_o$-connected that that the stabilizer of $(e^E)\tilde{}~$ 
is not larger than $M_K(E)$. In this case $\rho(M(E))(e^E)\tilde{}~$ 
is an imbedded image of $X^E$. Thus, for every $E_o$-connected $E$
we have an imbedding of $X^E$ into the boundary which we will denote 
by $\iota_E$. It is interesting to note that, in terms of the vector space
$V'$ (identified with its image in $PV$), the set $\iota_E(X^E)$ 
is just the parallel translate by $(e^E)\tilde{}~$ of $X^E = M(E)\tilde{e}$.
(Observe that $\tilde{e}$ is now identified with $0 \in V'$.) 

 It is easy to show that, for any $E_0$-connected
$E$, there is a unique maximal set $E'$ such that $E'$ is the
union of $E$ and of some components contained in $E_0$. Writing
$E'=E\cup E''$, we have $X^{E'}=X^E\times X^{E''}$. For all $H\in
\mathfrak{a}(E')^+$, when the image of $X^{E'}$ in $PV$ is
translated by $\rho(\exp tH)$ and we let $t\rightarrow \infty$,
the limit will be $\iota_E(X^E)$, while $X^{E''}$ will contract
and disappear.

We also note that applying $\rho(\exp tH)$ and letting $t\rightarrow
\infty$ actually moves all points of $X$ into $\iota_E(X^E)$. This
follows easily from $X=M(E')A(E')X^{E'}$ which, in turn, is a
consequence of the Iwasawa decomposition.

It is easy to see that the boundary consists of the $K$-images
of the sets $\iota_E(X^E)$. In fact if $\{k_\nu a_\nu\cdot o\}$ is
a sequence of points in $X$ tending to infinity, by compactness it
has a subsequence tending to a point in $k\cdot \iota_E(X^E)$ for
some $E$ and some $k\in K$.

To determine the stabilizers of boundary points, let $E$ be
$E_0$-connected and let $H\in \mathfrak{a}(E')^+$. Then the
stabilizer of $\exp tH\cdot o$ is $K^{\exp tH}$. As $t\rightarrow
\infty$, this group gets deformed into $M_K(E')N(E')$. (Indeed, 
any element of $M_K(E')N(E')$ can be written 
$mn=m \exp \sum_\gamma X^\gamma$, ($X^\gamma \in \mathfrak{g},   
(\gamma|H)>0)$ and is the limit of $k_t^{\exp tH}$ with 
$k_t=m\exp\sum e^{-t(\gamma|H)}(X^\gamma +\theta X^\gamma)$ in $K$.
As explained in [10, Sec. 3] this is the basic phenomenon behind
Bolyai's and Lobachevsky's definition of horicycles.) It follows
easily that the stabilizer of $\iota_E(o)$ is $M_K(E')A(E')N(E')$
and the stabilizer of  $\iota_E(X^E)$ is $B(E')$.

In this way all the Satake axioms ([13, p. 100], [10, Sec. 4]) are
verified, so our construction gives exactly the same
compactifications as the original one of Satake.

\bibliographystyle{plain}

\begin{flushright}
{Adam Kor\'{a}nyi\\ Mathematics Department\\ H. H. Lehman
College\\
Bronx, NY 10468, USA}
\end{flushright}

\end{document}